\documentclass[12pt]{article}
\usepackage{amsmath}
\usepackage{amsfonts}
\usepackage{geometry}

\input{tcilatex}
\begin{document}

\title{\textbf{Some Characterizations on the Normalized Lommel, Struve and
Bessel Functions of the First Kind}}
\author{Rabha M. El-Ashwah$^{\text{1}}$\ and Alaa H. El-Qadeem$^{\text{2}}$
\and $^{\text{1}}${\small Department of Mathematics, Faculty of Science,
Damietta University, New Damietta 34517, Egypt} \and $^{\text{2}}${\small %
Department of Mathematics, Faculty of Science, Zagazig University, Zagazig
44519, Egypt} \and {\small r\_elashwah@yahoo.com} \\
{\small alaahassan1986@yahoo.com \& alaahassan@zu.edu.eg}}
\date{}
\maketitle

\begin{abstract}
In this paper, we introduce new technique for determining some necessary and
sufficient conditions of the normalized Bessel functions $j_{\nu }$,
normalized Struve functions $h_{\nu }$ and normalized Lommel functions $%
s_{\mu ,\nu }$ of the first kind, to be in the subclasses of starlike and
convex functions of order $\alpha $ and type $\beta $.
\end{abstract}



\section*{1. Introduction}

Let $\mathcal{A}$ denotes the class of analytic functions of the form%
\begin{equation}
f(z)=z+\dsum\limits_{k=2}^{\infty }a_{k}z^{k},  \tag{1.1}
\end{equation}%
which are defined on the open unit disc $U=\{z\in
\mathbb{C}
:|z|$ $<1\}$ and satisfy the normalization conditions $f(0)=f^{\prime
}(0)-1=0$, and let $T$ be the subclass of $\mathcal{A}$ consisting of
functions of the form%
\begin{equation}
f(z)=z-\dsum\limits_{k=2}^{\infty }a_{k}z^{k}\text{ \ }(a_{k}\geq 0).
\tag{1.2}
\end{equation}%
\newline

\noindent \textbf{Definition 1.} \textit{For }$0\leq \alpha <1$\textit{\ and
}$0<\beta \leq 1,$\textit{\ let }$S^{\ast }(\alpha ,\beta )$\textit{\
denotes the subclass of }$\mathcal{A}$\textit{\ consisting of functions }$%
f(z)$\textit{\ of the form (1.1) and satisfy}%
\begin{equation}
\left\vert \frac{\frac{zf^{\prime }(z)}{f(z)}-1}{\frac{zf^{\prime }(z)}{f(z)}%
+\left( 1-2\alpha \right) }\right\vert <\beta ;\text{ }z\in U,  \tag{1.3}
\end{equation}%
\textit{and }$f\in K(\alpha ,\beta )$\textit{\ denotes the subclass of }$%
\mathcal{A}$\textit{\ consisting of functions }$f(z)$\textit{\ of the form
(1.1) such that}%
\begin{equation}
\left\vert \frac{\frac{zf^{\prime \prime }(z)}{f^{\prime }(z)}}{\frac{%
zf^{\prime \prime }(z)}{f^{\prime }(z)}+2\left( 1-\alpha \right) }%
\right\vert <\beta ;\text{ }z\in U.  \tag{1.4}
\end{equation}

\noindent The subclasses $S^{\ast }(\alpha ,\beta )$ and $K(\alpha ,\beta )$
are the well-known subclasses of starlike and convex functions of order $%
\alpha $ and type $\beta $, respectively, introduced by Gupta and Jain \cite%
{Gupta}.\newline
Moreover, let $T^{\ast }(\alpha ,\beta )$ and $C(\alpha ,\beta )$ be two
subclasses of $T$ defined by%
\begin{equation*}
T^{\ast }(\alpha ,\beta )=S^{\ast }(\alpha ,\beta )\cap T\ \text{and}\
C(\alpha ,\beta )=K(\alpha ,\beta )\cap T.
\end{equation*}%
\noindent It is known that:%
\begin{equation*}
f(z)\in C(\alpha ,\beta )\Longleftrightarrow zf^{\prime }(z)\in T^{\ast
}(\alpha ,\beta ),
\end{equation*}%
\bigskip

\noindent We note that $S^{\ast }(\alpha ,1)=S^{\ast }(\alpha )$ and $%
K(\alpha ,1)=K(\alpha )$ the subclasses of starlike and convex functions of
order $\alpha (0\leq \alpha <1),$ respectively, which were introduced by
Robertson \cite{Robertson}, Schild \cite{Schild} and MacGregor \cite%
{MacGregor}. Also, $T^{\ast }(\alpha ,1)=T^{\ast }(\alpha )$ and$\ C(\alpha
,1)=C(\alpha )$, are the subclasses of starlike and convex functions of
order $\alpha $ with negative coefficients introduced by see Silverman \cite%
{Silverman1}.

\bigskip

\noindent Recently, several researchers studied some subclasses of analytic
functions, $\tciFourier \subset \mathcal{A}$, involving special functions,
to find different conditions such that the members of $\tciFourier $ have
certain geometric properties such as univalency, starlikeness or convexity
in $U$. In this context many results are available in the literature
regarding the generalized hypergeometric functions (see \cite{Mocanu Miller}
and \cite{Singh Ruscheweyh}). Special functions, like Bessel, Struve and
Lommel functions of the first kind have many attractive geometric
properties. Recently, the geometric properties of these special functions
were investigated motivated by some earlier results. In the sixties Brown
\cite{Brown1}-\cite{Brown3}, Kreyszig and Todd \cite{Kreyszig Todd} and in
Wilf \cite{Wilf} the univalence and starlikeness of Bessel functions of the
first kind was considered, while in the recent years the radii of
univalence, starlikeness and convexity for the normalized forms of Bessel,
Struve and Lommel functions of the first kind were obtained, see the papers
\cite{Aktas}-\cite{Baricz Yagmur}, \cite{Szasz}, \cite{Szasz Kupan}, \cite%
{H1}-\cite{H4} and the references cited therein. In the above works the
authors used intensively some properties of the positive zeros of Bessel,
Struve and Lommel functions of the first kind, under some conditions. In
this paper our aim is to give some new results for the starlikeness and
convexity of the normalized Bessel, Struve and Lommel functions of the first
kind. \newline
In this paper, we consider three classical special functions, the Bessel
function of the first kind $J_{\nu }$, the Struve function of the first kind
$H_{\nu }$ and the Lommel function of the first kind $S_{\mu ,\nu }$. It is
known that the Bessel functions has the infinite series representation%
\begin{equation}
J_{\nu }(z)=\dsum\limits_{k=0}^{\infty }\frac{(-1)^{k}}{\Gamma \left(
k+1\right) \Gamma \left( k+\nu +1\right) }\left( \frac{z}{2}\right) ^{2k+\nu
},  \tag{1.5}
\end{equation}%
\noindent where $z,\nu \in
\mathbb{C}
$ such that $\nu \notin
\mathbb{Z}
^{-}:=\left\{ -1,-2,...\right\} $. Also, the Struve and Lommel functions can
be represented as the infinite series%
\begin{equation}
H_{\nu }(z)=\dsum\limits_{k=0}^{\infty }\frac{(-1)^{k}}{\Gamma \left( k+%
\frac{3}{2}\right) \Gamma \left( k+\nu +\frac{3}{2}\right) }\left( \frac{z}{2%
}\right) ^{2k+\nu +1};\text{ }\nu +\frac{1}{2}\notin
\mathbb{Z}
^{-},  \tag{1.6}
\end{equation}%
and%
\begin{equation}
S_{\mu ,\nu }(z)=\frac{z^{\mu +1}}{4}\dsum\limits_{k=0}^{\infty }\frac{%
(-1)^{k}\Gamma \left( \frac{\mu -\nu +1}{2}\right) \Gamma \left( \frac{\mu
+\nu +1}{2}\right) }{\Gamma \left( k+\frac{\mu -\nu +3}{2}\right) \Gamma
\left( k+\frac{\mu +\nu +3}{2}\right) }\left( \frac{z}{2}\right) ^{2k};\text{
}\frac{\mu \pm \nu +1}{2}\notin
\mathbb{Z}
^{-},  \tag{1.7}
\end{equation}%
for $\mu ,\nu ,z\in
\mathbb{C}
.$ In addition, we know that the Bessel function $J_{\nu }(z)$ is a solution
of the homogeneous Bessel differential equation (see \cite[p. 217]{Olver})%
\begin{equation*}
z^{2}w^{\prime \prime }(z)+zw^{\prime }(z)+\left( z^{2}-\nu ^{2}\right)
w(z)=0,
\end{equation*}

\noindent and the Struve function $H_{\nu }(z)$ and Lommel function\ $S_{\mu
,\nu }(z)$ are a particular solutions of the following non-homogeneous
Bessel differential equations (see \cite[p. 288 and p. 294]{Olver}),
respectively%
\begin{equation*}
z^{2}w^{\prime \prime }(z)+zw^{\prime }(z)+\left( z^{2}-\nu ^{2}\right) w(z)=%
\frac{\left( \frac{z}{2}\right) ^{\nu -1}}{\sqrt{\pi }\Gamma \left( \nu +%
\frac{1}{2}\right) },
\end{equation*}%
\noindent and%
\begin{equation*}
z^{2}w^{\prime \prime }(z)+zw^{\prime }(z)+\left( z^{2}-\nu ^{2}\right)
w(z)=z^{\mu +1}.
\end{equation*}%
It is worthy to notice that the functions $J_{\nu }(z),\,H_{\nu }(z)$ and $%
S_{\mu ,\nu }(z)$ are explicitly defined in terms of the hypergeometric
function $_{1}F_{2}$ by the following:%
\begin{equation*}
J_{\nu }(z)=\frac{\left( \frac{z}{2}\right) ^{\nu }}{\Gamma \left( \nu
+1\right) }%
\begin{array}{c}
_{1}F_{2}\left( 1;1,\nu +1;-\frac{z^{2}}{4}\right)%
\end{array}%
;\text{ }\nu \notin
\mathbb{Z}
^{-},
\end{equation*}%
\begin{equation*}
H_{\nu }(z)=\frac{\left( \frac{z}{2}\right) ^{\nu +1}}{\sqrt{\frac{\pi }{4}}%
\Gamma \left( \nu +\frac{3}{2}\right) }%
\begin{array}{c}
_{1}F_{2}\left( 1;\frac{3}{2},\nu +\frac{3}{2};-\frac{z^{2}}{4}\right)%
\end{array}%
;\text{ }\nu +\tfrac{1}{2}\notin
\mathbb{Z}
^{-},
\end{equation*}%
and%
\begin{equation*}
S_{\mu ,\nu }(z)=\frac{z^{\mu +1}}{\left( \mu -\nu +1\right) \left( \mu +\nu
+1\right) }%
\begin{array}{c}
_{1}F_{2}\left( 1;\frac{\mu -\nu +3}{2},\frac{\mu +\nu +3}{2};-\frac{z^{2}}{4%
}\right)%
\end{array}%
;\text{ }\tfrac{\mu \pm \nu +1}{2}\notin
\mathbb{Z}
^{-}.
\end{equation*}

\noindent We refer to Watson's treatise \cite{Watson} for comprehensive
information about these functions.

\bigskip

\noindent In this paper, we are mainly interested in the normalized Bessel
function of the first kind $j_{\nu }:U\rightarrow C$, normalized Struve
function of the first kind $h_{\nu }:U\rightarrow C$, and normalized Lommel
functions of the first kind $s_{\mu ,\nu }:U\rightarrow C$, which are
defined as follows%
\begin{equation}
j_{\nu }(z):=\Gamma \left( \nu +1\right) z^{1-\frac{\nu }{2}}J_{\nu }(2\sqrt{%
z})=\dsum\limits_{k=0}^{\infty }\frac{(-1)^{k}}{\left( 1\right) _{k}\left(
\nu +1\right) _{k}}z^{k+1};\text{ }\nu \notin
\mathbb{Z}
^{-},  \tag{1.8}
\end{equation}%
\begin{equation}
h_{\nu }(z):=\Gamma \left( \tfrac{3}{2}\right) \Gamma \left( \nu +\tfrac{3}{2%
}\right) z^{1-\frac{\nu }{2}}H_{\nu }(2\sqrt{z})=\dsum\limits_{k=0}^{\infty }%
\frac{(-1)^{k}}{\left( \tfrac{3}{2}\right) _{k}\left( \nu +\tfrac{3}{2}%
\right) _{k}}z^{k+1};\text{ }\nu +\tfrac{1}{2}\notin
\mathbb{Z}
^{-},  \tag{1.9}
\end{equation}%
and%
\begin{equation}
s_{\mu ,\nu }(z):=\left( \mu {\small -}\nu {\small +}1\right) \left( \mu
{\small +}\nu {\small +}1\right) z^{-\mu }S_{\mu ,\nu }(2\sqrt{z}%
)=\dsum\limits_{k=0}^{\infty }\frac{(-1)^{k}}{\left( \frac{\mu -\nu +3}{2}%
\right) _{k}\left( \frac{\mu +\nu +3}{2}\right) _{k}}z^{k+1};\text{ }\tfrac{%
\mu \pm \nu +1}{2}\notin
\mathbb{Z}
^{-}.  \tag{1.10}
\end{equation}%
Observe that%
\begin{equation}
s_{\nu ,\nu }(z)=h_{\nu }(z)\text{ \ and \ }s_{\nu -1,\nu }(z)=j_{\nu }(z).
\tag{1.11}
\end{equation}

\noindent Very recently, Cho et al. \cite{Cho Lee Srivastava} and
Murugusundaramoorthy and Janani \cite{Murugusundaramoorthy Janani} (see also
Porwal and Dixit \cite{Porwal Dixit}) introduced some characterization of
generalized Bessel functions of first kind to be in certain subclasses of
uniformly starlike and uniformly convex functions. In the present paper, we
determine necessary and sufficient conditions for normalized Bessel function
of the first kind, normalized Struve function of the first kind and
normalized Lommel functions of the first kind to be in certain Subclasses of
analytic functions.

\section*{\protect\large 2. Characterizations on Lommel Functions}

Unless otherwise mentioned, we assume in the reminder of this paper that$,$ $%
0\leq \alpha <1,$ $0<\beta \leq 1,$ and $z\in U$.\newline
To establish our main results, we shall require the following lemmas:

\noindent \textbf{Lemma 1} (\cite{Gupta}). \textit{(i) A sufficient
condition for a function }$f$\textit{\ of the form (1.1) to be in the class }%
$S^{\ast }(\alpha ,\beta )$\textit{\ is that}%
\begin{equation}
\dsum\limits_{k=2}^{\infty }\left[ k-1+\beta \left( k+1-2\alpha \right) %
\right] \left\vert a_{k}\right\vert \leq 2\beta \left( 1-\alpha \right) .
\tag{2.1}
\end{equation}%
\textit{(ii) A necessary and sufficient condition for a function }$f$\textit{%
\ of the form (1.2) to be in the }$T^{\ast }(\alpha ,\beta )$\textit{\ is
that}%
\begin{equation}
\dsum\limits_{k=2}^{\infty }\left[ k-1+\beta \left( k+1-2\alpha \right) %
\right] a_{k}\leq 2\beta \left( 1-\alpha \right) .  \tag{2.2}
\end{equation}%
\textit{\noindent }\textbf{Lemma 2 (\cite{Gupta}).}\textit{\ (i) A
sufficient condition for a function }$f$\textit{\ of the form (1.1) to be in
the class }$K(\alpha ,\beta )$\textit{\ is that}%
\begin{equation}
\dsum\limits_{k=2}^{\infty }k\left[ k-1+\beta \left( k+1-2\alpha \right) %
\right] \left\vert a_{k}\right\vert \leq 2\beta \left( 1-\alpha \right) .
\tag{2.3}
\end{equation}%
\textit{(ii) A necessary and sufficient condition for a function }$f$\textit{%
\ of the form (1.2) to be in the }$C(\alpha ,\beta )$\textit{\ is that}%
\begin{equation}
\dsum\limits_{k=2}^{\infty }k\left[ k-1+\beta \left( k+1-2\alpha \right) %
\right] a_{k}\leq 2\beta \left( 1-\alpha \right) .  \tag{2.4}
\end{equation}

\noindent If
\begin{equation}
f_{0}(z)=\frac{z}{1+z}=\tsum\nolimits_{k=0}^{\infty }\left( -1\right)
^{k}z^{k+1}\left( z\in U\right)  \tag{2.5}
\end{equation}%
and using the convolution principle, then let us define the function
\begin{equation*}
\mathfrak{s}_{\mu ,\nu }(z):=s_{\mu ,\nu }(z)\ast f_{0}(z),
\end{equation*}%
then we have the first result as follows.

\noindent \textbf{Theorem 1}. \textit{If }$\mu >\nu -3$\textit{, then the
condition}%
\begin{equation}
\left( 1+\beta \right) \mathfrak{s}_{\mu +2,\nu }^{\prime }(1)+2\beta \left(
1-\alpha \right) \mathfrak{s}_{\mu +2,\nu }(1)\leq \frac{8\beta (1-\alpha )}{%
\left( \mu -\nu +3\right) \left( \mu +\nu +3\right) }  \tag{2.6}
\end{equation}%
\textit{suffices that }$\mathfrak{s}_{\mu ,\nu }(z)\in S^{\ast }(\alpha
,\beta )$\textit{.}

\noindent \textit{Proof.} Since%
\begin{equation*}
\mathfrak{s}_{\mu ,\nu }(z)=z+\dsum\limits_{k=2}^{\infty }\frac{1}{\left(
\frac{\mu -\nu +3}{2}\right) _{k-1}\left( \frac{\mu +\nu +3}{2}\right) _{k-1}%
}z^{k};\text{ }\tfrac{\mu \pm \nu +1}{2}\notin
\mathbb{Z}
^{-}
\end{equation*}%
By virtue of (i) in Lemma 1, it is suffices to show that
\begin{equation*}
\dsum\limits_{k=2}^{\infty }\left[ k-1+\beta \left( k+1-2\alpha \right) %
\right] \frac{1}{\left( \frac{\mu -\nu +3}{2}\right) _{k-1}\left( \frac{\mu
+\nu +3}{2}\right) _{k-1}}\leq 2\beta (1-\alpha ).
\end{equation*}

\noindent By simple computation, we have%
\begin{eqnarray*}
\tciLaplace (\alpha ,\beta ;\mu ,\nu ) &=&\dsum\limits_{k=2}^{\infty }\left[
k-1+\beta \left( k+1-2\alpha \right) \right] \tfrac{1}{\left( \frac{\mu -\nu
+3}{2}\right) _{k-1}\left( \frac{\mu +\nu +3}{2}\right) _{k-1}} \\
&=&\dsum\limits_{k=0}^{\infty }\left[ k+1+\beta \left( k+3-2\alpha \right) %
\right] \tfrac{1}{\left( \frac{\mu -\nu +3}{2}\right) _{k+1}\left( \frac{\mu
+\nu +3}{2}\right) _{k+1}} \\
&=&\dsum\limits_{k=0}^{\infty }\left[ \left( k+1\right) \left( 1+\beta
\right) +2\beta \left( 1-\alpha \right) \right] \tfrac{1}{\left( \frac{\mu
-\nu +3}{2}\right) _{k+1}\left( \frac{\mu +\nu +3}{2}\right) _{k+1}} \\
&=&\left( 1+\beta \right) \dsum\limits_{k=0}^{\infty }\left( k+1\right)
\tfrac{1}{\left( \frac{\mu -\nu +3}{2}\right) _{k+1}\left( \frac{\mu +\nu +3%
}{2}\right) _{k+1}}+2\beta \left( 1-\alpha \right)
\dsum\limits_{k=0}^{\infty }\tfrac{1}{\left( \frac{\mu -\nu +3}{2}\right)
_{k+1}\left( \frac{\mu +\nu +3}{2}\right) _{k+1}}
\end{eqnarray*}%
\begin{eqnarray}
&=&\tfrac{\left( 1+\beta \right) }{\left( \frac{\mu -\nu +3}{2}\right)
\left( \frac{\mu +\nu +3}{2}\right) }\dsum\limits_{k=0}^{\infty }\left(
k+1\right) \tfrac{1}{\left( \frac{\mu -\nu +5}{2}\right) _{k}\left( \frac{%
\mu +\nu +5}{2}\right) _{k}}+\tfrac{2\beta \left( 1-\alpha \right) }{\left(
\frac{\mu -\nu +3}{2}\right) \left( \frac{\mu +\nu +3}{2}\right) }%
\dsum\limits_{k=0}^{\infty }\tfrac{1}{\left( \frac{\mu -\nu +5}{2}\right)
_{k}\left( \frac{\mu +\nu +5}{2}\right) _{k}}  \notag \\
&=&\tfrac{\left( 1+\beta \right) }{\left( \frac{\mu -\nu +3}{2}\right)
\left( \frac{\mu +\nu +3}{2}\right) }\dsum\limits_{k=0}^{\infty }\left(
k+1\right) \tfrac{1}{\left( \frac{\mu +2-\nu +3}{2}\right) _{k}\left( \frac{%
\mu +2+\nu +3}{2}\right) _{k}}+\tfrac{2\beta \left( 1-\alpha \right) }{%
\left( \frac{\mu -\nu +3}{2}\right) \left( \frac{\mu +\nu +3}{2}\right) }%
\dsum\limits_{k=0}^{\infty }\tfrac{1}{\left( \frac{\mu +2-\nu +3}{2}\right)
_{k}\left( \frac{\mu +2+\nu +3}{2}\right) _{k}}  \notag \\
&=&\tfrac{\left( 1+\beta \right) }{\left( \frac{\mu -\nu +3}{2}\right)
\left( \frac{\mu +\nu +3}{2}\right) }\dsum\limits_{k=0}^{\infty }\left(
k+1\right) \tfrac{1}{\left( \frac{\mu +2-\nu +3}{2}\right) _{k}\left( \frac{%
\mu +2+\nu +3}{2}\right) _{k}}+\tfrac{2\beta \left( 1-\alpha \right) }{%
\left( \frac{\mu -\nu +3}{2}\right) \left( \frac{\mu +\nu +3}{2}\right) }%
\dsum\limits_{k=0}^{\infty }\tfrac{1}{\left( \frac{\mu +2-\nu +3}{2}\right)
_{k}\left( \frac{\mu +2+\nu +3}{2}\right) _{k}}  \notag \\
&=&\left( \tfrac{\mu -\nu +3}{2}\right) \left( \tfrac{\mu +\nu +3}{2}\right) %
\left[ \left( 1+\beta \right) \mathfrak{s}_{\mu +2,\nu }^{\prime }(1)+2\beta
\left( 1-\alpha \right) \mathfrak{s}_{\mu +2,\nu }(1)\underset{}{}\right] .
\TCItag{2.7}
\end{eqnarray}

\noindent Therefore, we see that the last expression (2.7) is bounded above
by $2\beta (1-\alpha )$ if condition (2.6) is satisfied. This completes the
proof of Theorem 1.

\bigskip

\noindent If
\begin{equation}
f_{1}(z)=z\left( 2-\frac{1}{1+z}\right) =z+\tsum\nolimits_{k=1}^{\infty
}\left( -1\right) ^{k+1}z^{k+1}\text{ \ }\left( z\in U\right) ,  \tag{2.8}
\end{equation}%
and the function
\begin{equation*}
\mathfrak{t}_{\mu ,\nu }(z):=s_{\mu ,\nu }(z)\ast f_{1}(z),
\end{equation*}%
then we have the following result.

\noindent \textbf{Theorem 2}. \textit{For }$\mu >\nu -3$\textit{. then the
condition (2.6) is the necessary and sufficient condition for }$t_{\mu ,\nu
}(z)$\textit{\ to be in the class }$T^{\ast }(\alpha ,\beta )$\textit{.%
\newline
Proof.} Since
\begin{equation}
\mathfrak{t}_{\mu ,\nu }(z)=z-\dsum\limits_{k=2}^{\infty }\frac{1}{\left(
\frac{\mu -\nu +3}{2}\right) _{k-1}\left( \frac{\mu +\nu +3}{2}\right) _{k-1}%
}z^{k};\text{ }\tfrac{\mu \pm \nu +1}{2}\notin
\mathbb{Z}
^{-},  \tag{2.9}
\end{equation}%
then by using Lemma 1 together with the same techniques given in the proof
of Theorem 1 , we have immediately Theorem 2.

\noindent \textbf{Theorem 3}. \textit{If }$\mu >\nu -3$\textit{, then the
condition}%
\begin{equation}
\left( 1+\beta \right) \mathfrak{s}_{\mu +2,\nu }^{^{\prime \prime
}}(1)+2\left( 1+\beta \right) \mathfrak{s}_{\mu +2,\nu }^{\prime
}(1)+2\left( 1-\alpha \right) \left( 2\beta -1\right) \mathfrak{s}_{\mu
+2,\nu }(1)\leq \frac{8\beta (1-\alpha )}{\left( \mu -\nu +3\right) \left(
\mu +\nu +3\right) }  \tag{2.10}
\end{equation}%
\textit{suffices that }$\mathfrak{s}_{\mu ,\nu }(z)\in K(\alpha ,\beta )$%
\textit{.}

\noindent \textit{Proof.} By virtue of Lemma 2, it is suffices to show that
\begin{equation*}
\dsum\limits_{k=2}^{\infty }k\left[ k-1+\beta \left( k+1-2\alpha \right) %
\right] \frac{1}{\left( \frac{\mu -\nu +3}{2}\right) _{k-1}\left( \frac{\mu
+\nu +3}{2}\right) _{k-1}}\leq 2\beta (1-\alpha ).
\end{equation*}

\noindent By simple computation, we have%
\begin{eqnarray}
&&\tciFourier (\alpha ,\beta ;\mu ,\nu )  \notag \\
&=&\dsum\limits_{k=2}^{\infty }k\left[ k-1+\beta \left( k+1-2\alpha \right) %
\right] \tfrac{1}{\left( \frac{\mu -\nu +3}{2}\right) _{k-1}\left( \frac{\mu
+\nu +3}{2}\right) _{k-1}}  \notag \\
&=&\dsum\limits_{k=0}^{\infty }\left( k+2\right) \left[ k+1+\beta \left(
k+1+2\left( 1-\alpha \right) \right) \right] \tfrac{1}{\left( \frac{\mu -\nu
+3}{2}\right) _{k+1}\left( \frac{\mu +\nu +3}{2}\right) _{k+1}}  \notag \\
&=&\dsum\limits_{k=0}^{\infty }\left[ \left( 1+\beta \right) \left(
k+1\right) \left( k\right) +2\left( 2+\beta -\alpha \right) \left(
k+1\right) +2\left( 1-\alpha \right) \left( 2\beta -1\right) \underset{}{}%
\right] \tfrac{1}{\left( \frac{\mu -\nu +3}{2}\right) _{k+1}\left( \frac{\mu
+\nu +3}{2}\right) _{k+1}}  \notag \\
&=&\left( 1+\beta \right) \dsum\limits_{k=0}^{\infty }\tfrac{\left(
k+1\right) \left( k\right) }{\left( \frac{\mu -\nu +3}{2}\right)
_{k+1}\left( \frac{\mu +\nu +3}{2}\right) _{k+1}}+2\left( 2+\beta -\alpha
\right) \dsum\limits_{k=0}^{\infty }\tfrac{k+1}{\left( \frac{\mu -\nu +3}{2}%
\right) _{k+1}\left( \frac{\mu +\nu +3}{2}\right) _{k+1}}  \notag \\
&&+2\left( 1-\alpha \right) \left( 2\beta -1\right)
\dsum\limits_{k=0}^{\infty }\tfrac{1}{\left( \frac{\mu -\nu +3}{2}\right)
_{k+1}\left( \frac{\mu +\nu +3}{2}\right) _{k+1}}  \notag \\
&=&\tfrac{1}{\left( \frac{\mu -\nu +3}{2}\right) \left( \frac{\mu +\nu +3}{2}%
\right) }\left[ \left( 1+\beta \right) \dsum\limits_{k=0}^{\infty }\tfrac{%
\left( k+1\right) \left( k\right) }{\left( \frac{\mu -\nu +5}{2}\right)
_{k}\left( \frac{\mu +\nu +5}{2}\right) _{k}}+2\left( 2+\beta -\alpha
\right) \dsum\limits_{k=0}^{\infty }\tfrac{k+1}{\left( \frac{\mu -\nu +5}{2}%
\right) _{k}\left( \frac{\mu +\nu +5}{2}\right) _{k}}\right.  \notag \\
&&\left. +2\left( 1-\alpha \right) \left( 2\beta -1\right)
\dsum\limits_{k=0}^{\infty }\tfrac{1}{\left( \frac{\mu -\nu +5}{2}\right)
_{k}\left( \frac{\mu +\nu +5}{2}\right) _{k}}\right]  \notag \\
&=&\tfrac{1}{\left( \frac{\mu -\nu +3}{2}\right) \left( \frac{\mu +\nu +3}{2}%
\right) }\left[ \left( 1{\small +}\beta \right) \mathfrak{s}_{\mu +2,\nu
}^{^{\prime \prime }}(1)+2\left( 1{\small +}\beta \right) \mathfrak{s}_{\mu
+2,\nu }^{\prime }(1)+2\left( 1{\small -}\alpha \right) \left( 2\beta
{\small -}1\right) \mathfrak{s}_{\mu +2,\nu }(1)\right]  \TCItag{2.11}
\end{eqnarray}

\noindent Therefore, we see that the expression (2.11) is bounded above by $%
2\beta (1-\alpha )$ if (2.10) is satisfied. Hence, the proof is completed.

\noindent \textbf{Theorem 4}. \textit{For }$\mu >\nu -3$\textit{. then the
condition (2.10) is the necessary and sufficient condition for }$t_{\mu ,\nu
}(z)$\textit{\ to be in the class }$C(\alpha ,\beta )$\textit{.}

\bigskip

\noindent Putting $\beta =1$ in Theorems 2, 4, we obtain the following
corollaries.

\noindent \textbf{Corollary 1.} \textit{The function }$t_{\mu ,\nu }(z)$%
\textit{\ is a starlike function of order }$\alpha \left( 0\leq \alpha
<1\right) $\textit{, if and only if}%
\begin{equation*}
\mathfrak{s}_{\mu +2,\nu }^{\prime }(1)+\left( 1-\alpha \right) \mathfrak{s}%
_{\mu +2,\nu }(1)\leq \frac{4(1-\alpha )}{\left( \mu -\nu +3\right) \left(
\mu +\nu +3\right) }.
\end{equation*}%
\noindent \textbf{Corollary 2.} \textit{The function }$t_{\mu ,\nu }(z)$%
\textit{\ is a convex function of order }$\alpha \left( 0\leq \alpha
<1\right) $\textit{, if and only if}%
\begin{equation*}
\mathfrak{s}_{\mu +2,\nu }^{^{\prime \prime }}(1)+2\mathfrak{s}_{\mu +2,\nu
}^{\prime }(1)+\left( 1-\alpha \right) \mathfrak{s}_{\mu +2,\nu }(1)\leq
\frac{4(1-\alpha )}{\left( \mu -\nu +3\right) \left( \mu +\nu +3\right) }.
\end{equation*}

\section*{\protect\large 3. Characterizations on Struve Functions}

\noindent Taking $\mu =\nu $ in Theorems 1-4, then we obtain the
corresponding results of Struve function $h_{\nu }$, as following:

\noindent \textbf{Theorem 5}. \textit{If }$\nu >-\frac{1}{2}$\textit{, then
the condition}%
\begin{equation}
\left( 1+\beta \right) \mathfrak{h}_{\nu +2}^{\prime }(1)+2\beta \left(
1-\alpha \right) \mathfrak{h}_{\nu +2}(1)\leq \frac{8\beta (1-\alpha )}{%
3\left( 2\nu +3\right) }  \tag{3.1}
\end{equation}%
\textit{suffices that }$h_{\nu }(z)\in S^{\ast }(\alpha ,\beta )$\textit{,
where }%
\begin{equation}
\mathfrak{h}_{\nu }(z)%
\begin{array}{c}
:=%
\end{array}%
h_{\nu }(z)\ast f_{0}(z)=z+\dsum\limits_{k=2}^{\infty }\frac{1}{\left(
\tfrac{3}{2}\right) _{k-1}\left( \nu +\tfrac{3}{2}\right) _{k-1}}%
z^{k+1}\left( \nu >-\frac{1}{2}\right) .  \tag{3.2}
\end{equation}%
$.$

\noindent \textbf{Theorem 6}. \textit{If }$\nu >-\frac{3}{2}$\textit{, then
the condition (3.1) is the necessary and sufficient condition for }$\hbar
_{\nu }(z)$\textit{\ to be in the class }$T^{\ast }(\alpha ,\beta )$\textit{%
, where }%
\begin{equation}
\mathcal{\hbar }_{\nu }(z)%
\begin{array}{c}
:=%
\end{array}%
\left( h_{\nu }(z)\ast f_{1}(z)\right) =z-\dsum\limits_{k=2}^{\infty }\frac{1%
}{\left( \tfrac{3}{2}\right) _{k-1}\left( \nu +\tfrac{3}{2}\right) _{k-1}}%
z^{k+1}\left( \nu >-\frac{1}{2}\right) .  \tag{3.3}
\end{equation}

\noindent \textbf{Theorem 7}. \textit{If }$\nu >-\frac{3}{2}$\textit{, then
the condition}%
\begin{equation}
\left( 1+\beta \right) \mathfrak{h}_{\nu +2}^{\prime \prime }(1)+2\left(
1+\beta \right) \mathfrak{h}_{\nu +2}^{\prime }(1)+2\left( 1-\alpha \right)
\left( 2\beta -1\right) \mathfrak{h}_{\nu +2}(1)\leq \frac{8\beta (1-\alpha )%
}{3\left( 2\nu +3\right) }  \tag{3.2}
\end{equation}%
\textit{suffices that }$\mathfrak{h}_{\nu }(z)\in K(\alpha ,\beta ).$

\noindent \textbf{Theorem 8}. \textit{If }$\nu >-\frac{3}{2}$\textit{, then
the condition (3.2) is the necessary and sufficient condition for }$\hbar
_{\nu }(z)$\textit{\ to be in the class }$C(\alpha ,\beta )$\textit{.}

\noindent Putting $\beta =1$ in Theorems 6, 8, we have the following
corollaries.

\bigskip

\noindent \textbf{Corollary 3.} \textit{The function }$\hbar _{\nu }(z)$%
\textit{\ is a starlike function of order }$\alpha \left( 0\leq \alpha
<1\right) $\textit{, if and only if}%
\begin{equation*}
\mathfrak{h}_{\nu +2}^{\prime }(1)+\left( 1-\alpha \right) \mathfrak{h}_{\nu
+2}(1)\leq \frac{4(1-\alpha )}{3\left( 2\nu +3\right) }\left( 0\leq \alpha
<1,\nu >-\frac{3}{2}\right) .
\end{equation*}%
\textbf{Corollary 4.} \textit{The function }$\hbar _{\nu }(z)$\textit{\ is a
convex function of order }$\alpha \left( 0\leq \alpha <1\right) $\textit{,
if and only if}%
\begin{equation*}
\mathfrak{h}_{\nu +2}^{\prime \prime }(1)+2\mathfrak{h}_{\nu +2}^{\prime
}(1)+\left( 1-\alpha \right) \mathfrak{h}_{\nu +2}(1)\leq \frac{4(1-\alpha )%
}{3\left( 2\nu +3\right) }\left( 0\leq \alpha <1,\nu >-\frac{3}{2}\right) .
\end{equation*}

\section*{\protect\large 4. Characterizations on Bessel Functions}

\noindent Taking $\mu =\nu -1$ in Theorems 1-4, then we obtain the
corresponding results of Bessel function $j_{\nu }$, as following:

\noindent \textbf{Theorem 9}. \textit{If }$\nu >-1$\textit{, then the
condition}%
\begin{equation}
\left( 1+\beta \right) \mathfrak{j}_{\nu +1}^{\prime }(1)+2\beta \left(
1-\alpha \right) \mathfrak{j}_{\nu +1}(1)\leq \frac{2\beta (1-\alpha )}{\nu
+1},  \tag{4.1}
\end{equation}%
\textit{suffices that }$j_{\nu }(z)\in S^{\ast }(\alpha ,\beta )$\textit{,
where}
\begin{equation}
\mathfrak{j}_{\nu }(z)%
\begin{array}{c}
:=%
\end{array}%
j_{\nu }(z)\ast f_{0}(z)=z+\dsum\limits_{k=2}^{\infty }\frac{1}{\left(
1\right) _{k-1}\left( \nu +1\right) _{k-1}}z^{k}\left( \nu >-1\right) .
\tag{4.2}
\end{equation}

\noindent \textbf{Theorem 10}. \textit{If }$\nu >-1$\textit{, then the
condition (4.1) is the necessary and sufficient condition for }$\jmath _{\nu
}(z)$\textit{\ to be in the class }$T^{\ast }(\alpha ,\beta )$\textit{,
where }%
\begin{equation}
\jmath _{\nu }(z)%
\begin{array}{c}
:=%
\end{array}%
j_{\nu }(z)\ast f_{1}(z)=z-\dsum\limits_{k=2}^{\infty }\frac{1}{\left(
1\right) _{k-1}\left( \nu +1\right) _{k-1}}z^{k}\left( \nu >-1\right) .
\tag{4.3}
\end{equation}

\noindent \textbf{Theorem 11}. \textit{If }$\nu >-1$\textit{, then the
condition}%
\begin{equation}
\left( 1+\beta \right) \mathfrak{j}_{\nu +1}^{^{\prime \prime }}(1)+2\left(
1+\beta \right) \mathfrak{j}_{\nu +1}^{\prime }(1)+2\left( 1-\alpha \right)
\left( 2\beta -1\right) \mathfrak{j}_{\nu +1}(1)\leq \frac{2\beta (1-\alpha )%
}{\nu +1},  \tag{4.4}
\end{equation}%
\textit{suffices that }$\mathfrak{j}_{\nu }(z)\in K(\alpha ,\beta ).$

\noindent \textbf{Theorem 12}. \textit{If }$\nu >-1$\textit{, then the
condition (4.4) is the necessary and sufficient condition for }$\mathcal{%
\jmath }_{\nu }(z)$ to be in the class $C(\alpha ,\beta )$.

\bigskip

\noindent Putting $\beta =1$ in Theorems 10, 12, we get the following
corollaries.

\noindent \textbf{Corollary 5.} \textit{The function }$\jmath _{\nu }(z)$%
\textit{\ is a starlike function of order }$\alpha \left( 0\leq \alpha
<1\right) $\textit{, if and only if}%
\begin{equation*}
\mathfrak{j}_{\nu +1}^{\prime }(1)+\left( 1-\alpha \right) \mathfrak{j}_{\nu
+1}(1)\leq \frac{1-\alpha }{\nu +1}.
\end{equation*}%
\noindent \textbf{Corollary 6.} \textit{The function }$\jmath _{\nu }(z)$%
\textit{\ is a convex function of order }$\alpha \left( 0\leq \alpha
<1\right) $\textit{, if and only if}%
\begin{equation*}
\mathfrak{j}_{\nu +1}^{^{\prime \prime }}(1)+2\mathfrak{j}_{\nu +1}^{\prime
}(1)+\left( 1-\alpha \right) \mathfrak{j}_{\nu +1}(1)\leq \frac{1-\alpha }{%
\nu +1}.
\end{equation*}

\end{document}